\documentclass{IEEEtran}

\usepackage[noadjust]{cite}
\usepackage{amsmath,amssymb,amsfonts}
\usepackage{amsthm}
\usepackage{algorithmic}
\usepackage{graphicx}
\usepackage{textcomp}
\usepackage{mathrsfs}
\usepackage{xcolor}

\newtheorem{thm}{Theorem}
\newtheorem{lem}[thm]{Lemma}
\newtheorem{cor}[thm]{Corollary}
\newtheorem{prop}[thm]{Proposition}
\theoremstyle{remark}
\newtheorem{rem}[thm]{Remark}
\newtheorem{assum}{Assumption}

\theoremstyle{definition}
\newtheorem{defn}{Definition}

\def\bpf{\begin{IEEEproof}}
\def\epf{\end{IEEEproof}}
\providecommand{\norm}[1]{\lVert#1\rVert}

\def\BibTeX{{\rm B\kern-.05em{\sc i\kern-.025em b}\kern-.08em
    T\kern-.1667em\lower.7ex\hbox{E}\kern-.125emX}}

\allowdisplaybreaks[4]

\begin{document}
\title{Convergence of stochastic nonlinear systems and implications for stochastic Model Predictive Control}
\author{Diego Mu\~{n}oz-Carpintero and Mark Cannon 
\thanks{This work was supported in part by FONDECYT-CONICYT  Postdoctorado N\textsuperscript{o} 3170040.}
\thanks{D.~Mu\~{n}oz-Carpintero is with Institute of Engineering Sciences, Universidad de O´Higgins, Rancagua,
2841959, Chile, and Department of
Electrical Engineering, FCFM, University of Chile, Santiago 8370451, Chile (e-mail: diego.munoz@uoh.cl).}
\thanks{M.~Cannon is with Department of Engineering Science, University of
Oxford, OX1 3PJ U.K. (e-mail: mark.cannon@eng.ox.ac.uk).}}

\maketitle

\begin{abstract}
The stability of stochastic Model Predictive Control (MPC) subject to additive disturbances is often demonstrated in the literature by constructing Lyapunov-like inequalities that ensure closed-loop performance bounds and boundedness of the state, but tight ultimate bounds for the state and non-conservative performance bounds are typically not determined. 
In this work we use an input-to-state stability property to find conditions that imply convergence with probability~1 of a disturbed nonlinear system to 
a minimal robust positively invariant set.
We discuss implications for the convergence of the state and control laws of stochastic MPC formulations, and we prove convergence results for several existing stochastic MPC formulations for linear and nonlinear systems.
\end{abstract}


\section{Introduction}
\label{sec:introduction}

Stochastic Model Predictive Control (MPC) takes account of stochastic disturbances affecting a system using knowledge of their probability distributions or samples generated by an oracle~\cite{kouvaritakisandcannon2015,mesbah2016}. Stochastic predictions generated by a system model are used to evaluate a control objective, which is typically the expected value of a sum of stage costs, and probabilistic constraints, which allow constraint violations up to some specified maximum probability. The motivation for this approach is to avoid the conservatism of worst-case formulations of robust MPC and to account for stochastic disturbances in the optimization of predicted performance~\cite{mayne2014}.

Stability analyses of stochastic MPC can be divided into cases in which disturbances are either additive or multiplicative. For regulation problems with multiplicative disturbances, Lyapunov stability and convergence of the system state to the origin can often be guaranteed \cite{cannonetal2009d,bernardinibemporad2012,farinaetal2016}. However, convergence of stochastic MPC with additive disturbances is often either ignored \cite{matuskoborrelli2012,fagianokhammash2012,schilbachetal2014}, or else analysed using Lyapunov-like functions without identifying ultimate bounds on the state and asymptotic properties of the control law \cite{paulsonetal2015,lietal2018,santosetal2019}.
In~\cite{cannonetal2009a,cannonetal2009b,kouvaritakisetal2010} this difficulty is tackled by redefining the cost function and determining asymptotic bounds on the time-average stage cost, which converges to the value associated with the unconstrained optimal control law. In~\cite{chatterjeeandlygeros2015} a detailed analysis of the convergence and performance of stochastic MPC techniques is provided, considering typical stability notions of Markov chains and exploring the effect on stability and performance of the properties of the resulting Lyapunov-like inequalities (or geometric drift conditions). In particular, certain Lyapunov-like inequalities provide the robust notion of input-to-state stability (ISS), which may be used to derive (potentially conservative) ultimate boundedness conditions. 
These results are used for example in~\cite{paulsonetal2015,mishraetal2016} to analyse the stability of particular stochastic MPC formulations.

A difficulty with the stability analyses discussed so far is that stability and convergence properties of stochastic MPC may be stronger than what can be inferred directly from Lyapunov-like conditions. These may imply bounds on stored energy, measured for example by a quadratic functional of the system state~\cite{paulsonetal2015}, or that the state is bounded~\cite{goulartandkerrigan2008} without providing a tight limit set. Likewise, limit average performance bounds may be found~\cite{kouvaritakisetal2013} without guarantees that these bounds are tight.

However~\cite{lorenzenetal2017} presents a stochastic MPC strategy and a convergence analysis showing almost sure convergence of the system to the minimal robust positively invariant (RPI) set associated with a terminal control law of the MPC formulation. This implies a bound on time-average performance, and moreover the ultimate bounds for the state and average performance are tight.
The approach of~\cite{lorenzenetal2017} uses the statistical properties of the disturbance input to show that, for any given initial condition and terminal set (which can be any set that is positively invariant under the terminal control law), there exists a finite interval on which the  state of the closed-loop system will reach this terminal set with a prescribed probability.
The Borel-Cantelli lemma 
is then used to show
that the terminal set (and therefore the minimal RPI set associated with the terminal control law) is reached with probability 1 on an infinite time interval. 
It follows that average performance converges with probability 1 to that associated with the MPC law on the terminal set. These results are derived from the properties of a  particular MPC strategy that may not apply to other MPC formulations; namely that if a predicted control sequence is feasible for the MPC optimization at a sufficiently large number of consecutive sampling instants, then the state necessarily reaches the desired terminal set. 

Here we extend the analysis of \cite{lorenzenetal2017} to general nonlinear stochastic systems with ISS Lyapunov inequalities. We analyse convergence of the state to any invariant subset of the closed-loop system state space under the assumption that arbitrarily small disturbances have a non-vanishing probability. We thus derive tight ultimate bounds on the state and, for the case that the closed-loop dynamics are linear within the limit set, non-conservative bounds on time-average performance. The convergence analysis applies to stochastic MPC algorithms that ensure ISS, either for input- or state-constrained linear systems (for which the closed-loop dynamics are generally nonlinear), or for nonlinear systems. We apply the analysis to the stochastic MPC strategies of~\cite{kouvaritakisetal2013,goulartetal2006} for linear systems and to that of~\cite{santosetal2019} for nonlinear systems, deriving new convergence results for these controllers. Particularly, we find tight ultimate bounds that had not been proved before: for the state, for the controllers of \cite{kouvaritakisetal2013,santosetal2019}; and for the time-average performance for the controller of \cite{kouvaritakisetal2013}.

Our analysis uses arguments similar to those in~\cite{lorenzenetal2017}. Specifically, for any given initial condition and any invariant set containing the origin of state space, an ISS Lyapunov function is used to show that the probability of the state not entering this set converges to zero on an infinite time interval. It can be concluded that the state converges with probability~1 to the minimal RPI set for the system. This analysis is used to show that the closed-loop system state under a stochastic MPC law satisfying a suitable Lyapunov inequality converges with probability~1 to the minimal RPI set associated with the MPC law. This set and the limit average performance can be evaluated non-conservatively (to any desired precision) if the dynamics within the set are linear.

An analysis yielding similar results was recently published in~\cite{munozandcannon2019}. Using results given in~\cite{meynandtweedie2009} on the convergence of Markov chains,~\cite{munozandcannon2019} demonstrates convergence (in distribution) to a stationary terminal distribution supported in the minimal RPI of the system. This convergence was proved by assuming linearity of the dynamics in this terminal set imposing a controllability assumption on the system within this terminal set.
On the other hand, the analysis presented here uses the properties of ISS Lyapunov functions to provide stronger convergence results that hold with probability~1, and moreover these are demonstrated without the need for linearity or controllability assumptions.
Nevertheless, for the special case of linear dynamics, we show here that it is possible to compute the tight ultimate bounds and non-conservative limit average performance bounds.


The structure of the paper is as follows. Section 2 introduces the setting of the problem. Section 3 presents the stability analysis for general nonlinear systems under the specified assumptions. Section 4 discusses the implications for stochastic MPC by applying the analysis to derive convergence properties of several existing stochastic MPC formulations (for systems with linear and nonlinear dynamics, and MPC strategies with implicit and explicit predicted terminal controllers). 
Finally, Section 5 provides concluding remarks. 

\subsection{Basic definitions and notation}

The sets of non-negative integers and non-negative reals are denoted ${\mathbb N}$ and ${\mathbb R}_+$, and $\mathbb{N}_{[a,b]}$ is the sequence $\{a,a+1,\ldots,b\}$ and $\mathbb{N}_k = \mathbb{N}_{[0,k]}$. For a sequence $\{x_0,x_1,\ldots\}$,
$x_{j|k}$ for $j\in\mathbb{N}$ denotes the predicted value of $x_{k+j}$ made at time $k$. For sets $X,Y\subseteq{\mathbb R}^n$, the Minkowski sum is given by $X\oplus Y=\{x+y : x\in X,\, y\in Y \}$. The Minkowski sum of a sequence of sets $\{X_j , \, j\in\mathbb{N}_k\}$ is denoted $\bigoplus_{j=0}^k X_j$. For $X\subseteq{\mathbb R}^n$, ${\bf 1}_X(x)$ is the indicator function of $x\in X$. The open unit ball in ${\mathbb R}^n$ is ${\mathbb B}$. A continuous function $\phi:{\mathbb R}_+ \rightarrow {\mathbb R}_+$ is a ${\mathcal K}$-function if it is continuous, strictly increasing and $\phi(0)=0$, and it is a ${\mathcal K}_\infty$-function if it is a ${\mathcal K}$-function and $\phi(s)\rightarrow \infty$ as $s\rightarrow \infty$. A continuous function $\phi:{\mathbb R}_+\times{\mathbb R}_+\rightarrow {\mathbb R}_+$ is a  ${\mathcal K} {\mathcal L}$-function if $\phi(\cdot,t)$ is a ${\mathcal K}$-function for all $t\in{\mathbb R}_+$ and if $\phi(s,\cdot)$ is decreasing for all $s\in{\mathbb R}_+$ with $\phi(s,t)\rightarrow 0$ as $t\rightarrow \infty$. The probability of an event $A$ is denoted ${\mathbb P}(A)$.
We use the term \textit{limit set} to refer to an invariant subset of state-space to which the system converges. 

\section{Problem Setting}\label{sec2}

Consider a discrete time nonlinear system given by
\begin{equation} \label{eq1} 
x_{k+1} =f(x_k,w_k),
\end{equation} 
where $x_{k}\in \mathbb{X}\subseteq{\mathbb R}^{n} $ is the state, $w_{k}\in\mathbb{W}\subseteq{\mathbb R}^{n_w}$ is the disturbance input, $f:\mathbb{X}\times\mathbb{W}\rightarrow \mathbb{X}$ is a function with ${f(0,0)=0}$, and $\mathbb{X}$ contains the origin in its interior. Current and future values of $w_k$ are unknown.

\begin{assum}\label{distass}
The disturbance sequence $\{w_0, w_1,\ldots\}$ is independent and identically distributed (i.i.d.) with ${\mathbb E}\{w_k\}=0$. The probability density function (PDF) of $w$ is supported in $\mathbb{W}$, a bounded set that contains the origin in its interior. Additionally, ${\mathbb P}\{\norm{w}\le \lambda \}>0$ for all $\lambda>0$.
\end{assum}

The requirement that disturbance inputs are i.i.d.\ and zero-mean is a standard assumption in stochastic MPC formulations in the literature. The assumption that ${\mathbb P}\{\norm{w}\le \lambda \}>0$ for all $\lambda>0$ clearly excludes certain disturbance distributions, but we note that it does not require continuity of the PDF and is satisfied by uniform and Gaussian distributions (or truncated Gaussian, to comply with Assumption \ref{distass}), among many others.



While the dynamics of system \eqref{eq1} are defined on $\mathbb{X}$, we consider conditions for convergence to 
an invariant set $\Omega\subseteq\mathbb{X}$, characterized in Assumption \ref{terminalass}.
The analysis of convergence is based on the notion of input-to-state stability~\cite{jiangandwang2001}.

%
 
\begin{assum}\label{terminalass}
The set $\Omega$ is bounded, contains the origin in its interior and is RPI, i.e.\  $f(x,w)\in\Omega$ for all $(x,w)\in\Omega\times\mathbb{W}$.
\end{assum}

\begin{defn}[Input-to-state stability] 
System \eqref{eq1} is input-to-state stable (ISS) with region of attraction $\mathbb{X}$
if there exist a ${\mathcal K}{\mathcal L}$-function $\beta(\cdot,\cdot)$ and a ${\mathcal K}$-function $\gamma(\cdot)$ such that, for all $k\in{\mathbb N}$, all $x_0\in\mathbb{X}$ and all admissible disturbance sequences $\{w_0,\ldots,w_{k-1}\} \in \mathbb{W}^k$, the state of~\eqref{eq1} satisfies $x_k\in\mathbb{X}$ and
\begin{equation} \label{iss}
\norm{x_k}\le \beta(\norm{x_0},k)+\gamma\left( \max_{j\in \{0,\ldots,k-1 \}} \norm{w_j}\right).
\end{equation}
\end{defn}

The following theorem provides a necessary and sufficient condition under which system (\ref{eq1}) is ISS.

\begin{thm}[\cite{jiangandwang2001}]\label{lemiss}
System \eqref{eq1} is ISS with region of attraction~$\mathbb{X}$  if, and only if, there exist ${\mathcal K}_\infty$-functions $\alpha_1(\cdot)$, $\alpha_2(\cdot)$, $\alpha_3(\cdot)$, a ${\mathcal K}$-function $\sigma(\cdot)$ and a continuous function $V:\mathbb{X}\rightarrow {\mathbb R}_+$ such that
\begin{subequations}
\begin{align}
& \alpha_1(\norm{x})\le V(x) \le \alpha_2(\norm{x}), \quad \forall x\in {\mathbb R}^{n},
\label{lyapkappa}
\\
& V\left(f(x,w)\right)-V(x)\le -\alpha_3(\norm{x})+\sigma(\norm{w}), 
\label{isslyap}
\end{align}
\end{subequations}
for all $(x,w)\in\mathbb{X}\times\mathbb{W}$.
In this case we say that $V(\cdot)$ is an ISS-Lyapunov function.
\end{thm}

\begin{assum}\label{lyapunov}
System \eqref{eq1} is ISS with region of attraction~$\mathbb{X}$.
\end{assum}

Assumption~\ref{lyapunov} implies that: (i) the origin of the state-space of the system $x_{k+1}=f(x_k,0)$ is asymptotically stable; (ii)  all trajectories of \eqref{eq1} are bounded since $\mathbb{W}$ is a bounded set; and (iii) 
all trajectories of \eqref{eq1} converge to the origin if $w_k\rightarrow0$ as $k\rightarrow\infty$. For details we refer the reader to~\cite{jiangandwang2001}.

We assume throughout this work that the system~\eqref{eq1} possesses an ISS-Lyapunov function, and hence that Assumption~\ref{lyapunov} holds. 
This assumption does not directly guarantee convergence to the set $\Omega$, however. 
Instead we combine this property with the stochastic nature of disturbances satisfying Assumption~\ref{distass}
to prove convergence to $\Omega$ in Section~\ref{sec3}. 

\section{Main result}\label{sec3}

As discussed above, ISS is not generally sufficient to determine non-conservative ultimate bounds for the state of~\eqref{eq1}. However, Section \ref{ssmainXf} shows that  under Assumptions~\ref{distass} and~\ref{lyapunov} the state of~\eqref{eq1} converges with probability 1 to any set $\Omega$ satisfying Assumption~\ref{terminalass}, and thus to the minimal RPI set defined by the intersection of all such sets.
%
Section \ref{secterminallinear} considers the special case of linear dynamics on $\Omega$ (for which the minimal RPI set can be determined with arbitrary precision); while~\cite{lorenzenetal2017} proves this for a particular MPC algorithm, in this section we extend the treatment to more general systems.

\subsection{Convergence to $\Omega$}\label{ssmainXf}

This section uses the ISS property of Assumption~\ref{lyapunov} to demonstrate almost sure convergence to a set $\Omega$ satisfying Assumption~\ref{terminalass} if the disturbance  satisfies Assumption~\ref{distass}.
The general idea is as follows. We note that for any state $x \notin \Omega$ there exists a set $W\subseteq \mathbb{W}$ such that the ISS-Lyapunov function decreases if $w\in W$ and such that $w\in W$ occurs with non-zero probability. It follows that, for any given $x_0\in \mathbb{X}$, there exists $N_f\in\mathbb{N}$ such that if $w\in W$ for $N_f$ consecutive time steps, then the Lyapunov function decreases enough to ensure that the state enters $\Omega$. This observation is used to show that, for any given $p\in(0,1]$, there exists a finite horizon over which $x$ reaches $\Omega$ with probability $1-p$. 
Finally, the Borell-Cantelli lemma is used to conclude that the state converges to $\Omega$ with probability 1, which implies that, for any $x_0\in \mathbb{X}$, there is zero probability of a disturbance sequence realisation $\{w_0,w_1,\ldots\}$ such that 
$x_k$ does not converge to $\Omega$ as $k\to\infty$. 

\begin{prop}\label{propW}
Under Assumption~\ref{lyapunov}, for given $\lambda \in (0,1)$ let $W(z) = \{w:\sigma(\norm{w})\le \lambda \alpha_3(z)\}$, where $\sigma(\cdot)$ and $\alpha_3(\cdot)$ are a $\mathcal{K}$-function and a $\mathcal{K}_\infty$-function satisfying \eqref{isslyap}, and where $z\in \mathbb{R}_+$. Then there exists $\lambda \in (0,1)$ and a ${\mathcal K}$-function $\xi(\cdot)$ such that the ISS-Lyapunov function $V(\cdot)$ satisfies
\begin{equation}\label{lyapstrict}
V\bigl( f(x,w) \bigr)- V(x) \leq -\xi(\norm{x}) ,
\end{equation}
whenever $w\in W(\norm{x})$.
\end{prop}
\bpf Note that $-\alpha_3(\norm{x})+\sigma(\norm{w})\le -(1-\lambda)\alpha_3(\norm{x})$ if $w\in W(\norm{x})$. It follows immediately that \eqref{isslyap} implies \eqref{lyapstrict} with $\xi(\norm{x})=(1-\lambda)\alpha_3(\norm{x})$. \epf
\begin{prop}\label{propNfXf}
Under Assumptions 
\ref{terminalass} and \ref{lyapunov}, for any $x_0\in \mathbb{X}$ there exists an integer $N_f=\lceil \alpha_2(r(x_0))/\xi(\epsilon) \rceil$ such that, if $w_j\in W(\epsilon)$ for $j=k,k+1,\ldots k+N_f-1$,
where $\epsilon = \sup \{ \rho :  \rho {\mathbb B}\subseteq \Omega\}$ and $k\in\mathbb{N}$ is arbitrary, then $x_{k+N_f}\in \Omega$.
\end{prop}

\bpf 
The positive invariance of $\Omega$ in Assumption~\ref{terminalass} implies that $x_{k+N}\in \Omega$ for all $N\in\mathbb{N}$ if $x_{k} \in \Omega$. Suppose therefore that $x_j\notin \Omega$ for $j=k,\ldots,k+N-1$ for given $N\in \mathbb{N}$. 
In this case $W(\norm{x_j}) \supseteq W(\epsilon)$. Then, if $w_j\in W(\epsilon)$ for $j=k,\ldots,k+N-1$, Proposition \ref{propW} implies
\[
\sum_{j=k}^{k+N-1} \xi(\norm{x_j}) \leq V(x_k) - V(x_{k+N}) ,
\]
and since Assumption~\ref{lyapunov} implies $V(x_k) \leq \alpha_2(\norm{x_k}) \leq \alpha_2(r(x_0))$, where $r(x_0) = \beta(\norm{x_0},0) + \gamma(\sup_{w\in\mathbb{W}}\norm{w})$ for some ${\mathcal K}{\mathcal L}$-function $\beta(\cdot,\cdot)$ and some ${\mathcal K}$-function $\gamma(\cdot)$, the right-hand side of this inequality can be no greater than $\alpha_2(r(x_0))$.
Furthermore, $x_j\notin \Omega$ implies $\xi(\norm{x_j}) \geq \xi(\epsilon)$ and hence $N \leq \alpha_2(r(x_0))/\xi(\epsilon)$. Choosing $N=N_f=\lceil \alpha_2(r(x_0))/\xi(\epsilon) \rceil$ therefore ensures that $x_{k+N_f}\in\Omega$. 
\epf
%


\begin{lem}\label{lemkey}
Under Assumptions \ref{distass}-\ref{lyapunov}, for any $x_0\in \mathbb{X}$ and any given $p\in (0,1]$ there exists an integer $N_p$ such that 
\begin{equation}\label{lemkeyeq}
{\mathbb P} \{ x_{N_p} \in \Omega \} \geq 1 - p.
\end{equation}
\end{lem}

\bpf
Proposition~\ref{propNfXf} and the positive invariance of $\Omega$ in Assumption~\ref{terminalass} ensure that, for any $x_0$ and $k\in{\mathbb N}$, ${x_{k+{N_f}}\in\Omega}$  whenever $w_j\in W(\epsilon)$ for all $j \in\{ k,\ldots,k+N_f-1\}$ with $N_f=\lceil \alpha_2(r(x_0))/\xi(\epsilon) \rceil$.
Let $p_\epsilon$ be the probability that $w_j\in W(\epsilon)$, then the i.i.d.\ property of Assumption~\ref{distass} implies $\mathbb{P}\{{ w_j \in W(\epsilon)}, \, j = k,\ldots,k+N_f-1\} = p_{\epsilon}^{N_f} > 0$ and, since the event that $w_j \notin W(\epsilon)$ for some $j \in \{k,\ldots,k+N_f-1\}$ does not guarantee that $x_{k+N_f}\notin\Omega$, we obtain the bound
\[
\mathbb{P} \{ x_{k+N_f}\!\notin\!\Omega\} \leq 1 - \mathbb{P} \bigl\{ w_j \!\in\! W(\epsilon), \, j = 1,\ldots,N_f\bigr\}  = 1 - p_{\epsilon}^{N_f} 
\]
Now consider $k > N_f$, and choose $j_0,j_1,\ldots , j_{\lfloor k/N_f \rfloor}$ so that $j_0 = 0$, $j_{\lfloor k/N_f \rfloor} \leq k$ and $j_{i+1} - j_i \geq N_f$ for all $i$. Then $x_{k}\notin\Omega$ only if $x_{j_i}\notin \Omega$ for all $i=0,\ldots, \lfloor N/N_f \rfloor$. Hence
\begin{equation}\label{eq:pAc}
\mathbb{P} \{ x_{k}\notin\Omega\} \leq (1 - p_\epsilon^{N_f})^{\lfloor k/ N_f \rfloor}  
\end{equation}
and (\ref{lemkeyeq}) therefore holds with $N_p = N_f \lceil \log p /\log (1-p_\epsilon^{N_f})\rceil$ if $p_\epsilon < 1$ or $N_p = N_f$ if $p_\epsilon = 1$.
\epf

Proposition~\ref{propNfXf} establishes that the state of (\ref{eq1}) reaches $\Omega$ in finite time if the disturbance input is small enough for a sufficiently large number of consecutive time steps. Lemma~\ref{lemkey} uses a lower bound on the probability of this event to show that the probability of the state reaching this set is arbitrarily close to 1 on a long enough (but finite) horizon. Furthermore, for any $k \in\mathbb{N}$, the argument used to prove Lemma~\ref{lemkey} implies
\[
\mathbb{P} \{ x_k \in \Omega \}  \geq 1 - (1 - p_\epsilon^{N_f})^{\lfloor k/ N_f \rfloor} ,
\]
and an immediate consequence is that $x_k$ converges to $\Omega$ as
\[
\lim_{k\to\infty} \mathbb{P} \{ x_k \in \Omega \}  = 1.
\]
We next use the Borel-Cantelli lemma (in a similar fashion to~\cite{lorenzenetal2017}) to prove the slightly stronger property that the state converges to $\Omega$ with probability 1.

\begin{thm}\label{mainthm1}
Under Assumptions \ref{distass}-\ref{lyapunov}, for any $x_0\in\mathbb{X}$ we have
\begin{equation}\label{eqthm1}
\mathbb{P} \Bigl\{ \lim_{k\rightarrow\infty}  1_{ \Omega}(x_k)=1 \Bigr\} =1.
\end{equation}
\end{thm}

\bpf
Let $A_k$ denote the event $x_k\notin \Omega$, then Lemma \ref{lemkey} implies (according to \eqref{eq:pAc})
\[
\sum_{k=0}^\infty \mathbb{P} \{ A_k\} \leq \sum_{k=0}^\infty \bigl( 1 - p_\epsilon^{N_f}\bigr)^{\lfloor k/N_f\rfloor} = N_fp_\epsilon^{-N_f} < \infty
\]
and the Borel-Cantelli lemma therefore implies 
\[
\mathbb{P} \biggl\{ \bigcap_{k=0}^\infty \bigcup_{j=k}^\infty A_j \biggr\} = 0 .
\]
But $A_{k+1} \subseteq A_k$ for all $k\in\mathbb{N}$ since $\Omega$ is RPI due to Assumption \ref{terminalass}. It follows that
$\mathbb{P} \{ \cap_{k=0}^\infty A_k \} = 0$,
which is equivalent to $\mathbb{P} \{ \lim_{k\to \infty} A_k \} = 0$ and hence~(\ref{eqthm1}).
\epf

\begin{defn}[Minimal RPI set] 
The minimal RPI set for system (\ref{eq1}) containing the origin, denoted $\mathbb{X}_\infty$, is defined as the intersection of all sets $X$ such that $0\in X \subseteq \mathbb{X}$ and $f(x,w)\in X$ for all $(x,w)\in X\times \mathbb{W}$.
\end{defn}

\begin{rem}~\label{rem:mRPI}
Theorem \ref{mainthm1} applies to any RPI set $\Omega$ that contains the origin in its interior.
In particular $\Omega=\mathbb{X}_\infty$ is the smallest set that can satisfy Assumption \ref{terminalass}~\cite[{Prop.~6.13}]{blanchini08}. And, since the state can escape any smaller set because it would not be invariant, $\mathbb{X}_\infty$ is the smallest set to which the state of (\ref{eq1}) converges with probability 1.
\end{rem}

\begin{cor}\label{mainthm2}
Let Assumptions~\ref{distass}-\ref{lyapunov} hold and let $x_0\in\mathbb{X}$. Then
\begin{equation}\label{convlimit}
\mathbb{P} \Bigl\{  \lim_{k\rightarrow\infty}  1_{ \mathbb{X}_\infty}(x_k)=1  \Bigr\} = 1 .
\end{equation}
\end{cor}
\bpf 
If Assumption \ref{terminalass} holds for any RPI set, then the minimal RPI set $\mathbb{X}_\infty$ exists and also satisfies this assumption. The result then follows directly from Theorem \ref{mainthm1}.
\epf

In this section we have demonstrated convergence with probability 1 of the state of~\eqref{eq1} to any RPI set containing the origin in its interior. 
Remark~\ref{rem:mRPI} and Corollary~\ref{mainthm2} therefore imply that the minimal RPI set $\mathbb{X}_\infty$ is a tight limit set of~\eqref{eq1}.
This improves on the result of~\cite{munozandcannon2019}, where convergence to $\mathbb{X}_\infty$ in probability is shown for the case that $f(x,w)=Ax+Dw$ for all $(x,w)\in\Omega\times\mathbb{W}$, where $(A,D)$ is controllable. 

\subsection{Convergence to a limit set with linear dynamics}
\label{secterminallinear}

Of particular interest when analysing stochastic MPC algorithms for constrained linear systems is the case in which the dynamics of system (\ref{eq1}) are linear on an RPI set containing the origin.
In this case the minimal RPI set $\mathbb{X}_\infty$ defining ultimate bounds for the state and the limit average performance can be explicitly determined.

\begin{assum}\label{terminallinearass}
There exists an RPI set $\Gamma\subseteq\mathbb{X}$, such that $f(x,w)=\Phi x+Dw$ for all $(x,w)\in\Gamma\times\mathbb{W}$, 
where $\Phi,D$  are matrices with  appropriate dimensions and $\Phi$ is Schur stable.
\end{assum}

Corollary~\ref{mainthm2} implies that the state of (\ref{eq1}) converges with probability~1 from any initial condition in $\Gamma$ to the minimal RPI set $\mathbb{X}_\infty$ given by
\begin{equation}\label{eq:mRPIset}
\mathbb{X}_\infty=\lim_{k\rightarrow\infty}\bigoplus_{j=1}^k \Phi^j D\mathbb{W} .
\end{equation}
This set can be computed with arbitrary precision (see e.g.~\cite{blanchini08} for details).

We next consider the asymptotic time-average value of a quadratic function of the state of (\ref{eq1}), representing, for example, a  quadratic performance cost. 


\begin{thm}\label{mainthm3}
Let Assumptions~\ref{distass}-\ref{terminallinearass} hold and let $x_0\in\mathbb{X}$. Then
\begin{equation}\label{perflimit}
\lim_{k\rightarrow\infty} \frac{1}{k}\sum_{j=0}^{k-1}{\mathbb E}\{x_j^\top Sx_j \} \leq l_{ss} 
\end{equation}
for any given $S=S^\top\succ 0$ where $l_{ss}={\mathbb E}\{w^\top D^\top PDw \}$ 
with $P\succ 0$ satisfying $P-\Phi^\top P\Phi = S$. 
\end{thm}

\bpf Let $V(x)=x^\top Px$. If $x_j\in\Gamma$, then 
\[
{\mathbb E}\{V(x_{j+1}) \}-{\mathbb E}\{V(x_j) \}= 
-{\mathbb E}\{x_j^\top Sx_j \}+{\mathbb E}\{w_j^\top D^\top PDw_j \}.
\]
On the other hand, if $x_j\notin \Gamma$, then
\begin{align*}
&{\mathbb E}\{V(x_{j+1}) \}-{\mathbb E}\{V(x_j) \} 
\\
&= {\mathbb E}\{f(x_j,w_j)^\top Pf(x_j,w_j)  \}-{\mathbb E}\{x_j^\top Px_j \} 
\\
& ={\mathbb E}\{\delta(x_j,w_j)\}\!+\!{\mathbb E}\{(\Phi x_{j\!}\!+\! D w_{j\!})^{\!\!\top\!} \! P (\Phi x_{j\!}\!+\! D w_{j\!} )\} \!-\!{\mathbb E}\{x_{j\!}^{\!\top} \! Px_{j\!}\} 
\\
& ={\mathbb E}\{\delta(x_j,w_j)\}-{\mathbb E}\{x_j^\top Sx_j \}+{\mathbb E}\{w_j^\top D^\top PDw_j \},
\end{align*} 
where $\delta(x_{j}, w_{j})  =  f(x_{j}, w_{j})^{\top} P f(x_{j}, w_{j}) - (\Phi x_{j}+ Dw_{j})^{\top} P(\Phi x_{j}+ Dw_{j})$ and $\mathbb{E}\{\delta(x_j,w_j)\}\le \nu$ for finite $\nu$ since $x_j$, $w_j$ are bounded due to \eqref{iss} and Assumption~\ref{distass}.
Therefore
\[
{\mathbb E}\{V(x_{j+1}) \}-{\mathbb E}\{V(x_j) \}\le 
\nu{\mathbb P}\{x_j\!\notin\! \Gamma\} - {\mathbb E}\{x_j^\top Sx_j \} + l_{ss}.
\]
Summing both sides of this inequality over $0\leq j < k$ yields
\[
{\mathbb E}\{V(x_{k}) \}-{\mathbb E}\{V(x_0) \} \!\leq\! 
\sum_{j=0}^{k-1}\bigl( \nu{\mathbb P}\{x_j\!\notin\! \Gamma \} -{\mathbb E}\{x_j^\top Sx_j \}\bigr) + kl_{ss} .
\]
Here
$\sum_{j=0}^{k-1}{\mathbb P}\{x_j\notin \Gamma \} \leq N_f p_\epsilon^{-N_f\!\!}$
and
${\mathbb E}\{V(x_{k}) \}$ is finite for all $k$ due to (\ref{iss}).
%
We therefore obtain (\ref{perflimit}) in the limit as $k\to\infty$. 
\epf






Under Assumptions \ref{distass}-\ref{terminallinearass},
the state of \eqref{eq1} therefore converges with probability 1 to 
the minimal RPI set $\mathbb{X}_\infty$ and the time-average performance converges to its limit average on this set.
Moreover the bound in \eqref{perflimit} is tight because $l_{ss}$ is equal to the time-average performance associated with the linear dynamics defined in Assumption \ref{terminallinearass}.

\section{Implications for stability and convergence of  stochastic MPC}\label{sec4}

This section uses the results of Section~\ref{sec3} to analyse the convergence of three existing stochastic MPC algorithms to a limit set of the closed loop system.
The first of these is for linear systems and assumes a control policy that is an affine function of the disturbance input~\cite{goulartetal2006,goulartandkerrigan2008}.
For this approach, convergence to a minimal RPI set was shown in~\cite{wangetal2008} by redefining the cost function and control policy; here we provide analogous results for the original formulation in~\cite{goulartandkerrigan2008}.
%
%
The second MPC algorithm is also for linear systems, but it assumes an affine disturbance feedback law with a different structure (striped and extending across an infinite prediction horizon), for which the gains are computed offline~\cite{kouvaritakisetal2013}. The proof of convergence for this second MPC formulation is provided here for the first time. The third is a generic stochastic MPC algorithm for nonlinear systems based on constraint-tightening~\cite{santosetal2019}. Only ISS is proved in~\cite{santosetal2019}; here we demonstrate convergence to a minimal RPI set.

The system dynamics for the nonlinear stochastic MPC formulation are defined in Section \ref{nlsmpc}.
On the other hand, Sections~\ref{sec:gandk} and~\ref{sec:kcandm} assume dynamics defined by 
\begin{equation}\label{lindyn}
x_{k+1}=Ax_k+Bu_k+Dw_k
\end{equation}
where $A,B,D$ are matrices of conformal dimensions, and $(A,B)$ is stabilizable. A measurement of the current state, $x_k$, is assumed to be available at time $k$, but current and future values of $w_k$ are unknown. In each case the disturbance sequence $\{w_0,w_1,\ldots\}$ is assumed to be i.i.d.\ with $\mathbb{E}\{w_k\} = 0$, and the PDF of $w$ is supported in a bounded set $\mathbb{W}$ containing the origin in its interior. These assumptions are included in~\cite{goulartandkerrigan2008,kouvaritakisetal2013, santosetal2019}; here we assume in addition that ${\mathbb P}\{\norm{w}\le \lambda \}>0$ for all $\lambda>0$ so that Assumption~\ref{distass} holds.


\subsection{Affine in the disturbance stochastic MPC}\label{sec:gandk}

In~\cite{goulartandkerrigan2008} the predicted control policy is an affine function of future disturbances. The expected value of a quadratic cost is minimized online subject to the condition that state and control constraints hold for all future realisations of disturbance inputs.
The state and control constraints take the form
\begin{equation}\label{goulart_cons}
(x_k,u_k)\in \mathbb{Z}
\end{equation}
for all $k\in\mathbb{N}$, where $\mathbb{Z}\subseteq {\mathbb R}^{n}\times{\mathbb R}^{n_u}$ is a convex and compact set containing the origin in its interior. 

The control input is determined at each discrete time instant by solving a stochastic optimal control problem. To avoid the computational load of optimizing an arbitrary feedback policy, predicted control inputs are parameterized for $i\in\mathbb{N}_{N-1}$ as
\begin{equation}
u_{i|k}=v_{i|k}+\sum_{j=0}^{i-1}M_{i,j|k}w_{j|k} ,
\end{equation}
where the open-loop control sequence
${\bf v}_k=\{v_{i|k}, i\in\mathbb{N}_{N-1}\}$ and feedback gains ${\bf M}_k =\{M_{i,j|k}, \, j\in\mathbb{N}_{i-1}, \, i\in\mathbb{N}_{[1,N-1]}\}$ are decision variables at time $k$.
%
For all $i\geq N$, predicted control inputs are defined  by $u_{i|k}=Kx_{i|k}$, where $A+BK$ is Schur stable.
The predicted cost at time $k$ is defined as
\[
J(x_k,{\bf v}_k,{\bf M}_k) ={\mathbb E}\Bigl\{x_{i|k}^\top Px_{i|k}+\sum_{i=0}^{N-1} 
(x_{i|k}^\top Qx_{i|k}+u_{i|k}^\top Ru_{i|k}) \Bigr\}
\]
where $Q \succeq 0$, $R\succ 0$, $(A,Q^{1/2})$ is detectable, 
and $P\succ 0$ is the solution of the algebraic Riccati equation $P=Q+A^\top PA-K^\top (R+B^\top PB)K$, $K=-(R+B^\top PB)^{-1}B^\top PA$~\cite{goulartandkerrigan2008}. 
%
A terminal constraint, $x_{N|k}\in\mathbb{X}_f$, is included in the optimal control problem,
where $\mathbb{X}_f$ is an RPI set for the system \eqref{lindyn} with control law $u_k=Kx_k$ and constraints $(x_k,Kx_k)\in{\mathbf Z}$. 

The optimal control problem solved at the $k$th instant is
\begin{align*}
\mathcal{P}_1: \ \ & \min_{{\bf v}_k,{\bf M}_k}
& & J(x_k,{\bf v}_k,{\bf M_k}) \quad \text{s.t.} \quad \forall w_{i|k}\in\mathbb{W}, \ \forall  i\in\mathbb{N}_{N-1} \\
&&& u_{i|k}=v_{i|k}+\sum_{j=0}^{i-1}M_{i,j|k}w_{j|k} \\
&&& (x_{i|k},u_{i|k})\in\mathbb{Z} \\
&&& x_{i+1|k}=Ax_{i|k}+Bu_{i|k}+Dw_{i|k}   \\
&&& x_{N|k}\in \mathbb{X}_f \\
& && x_{0|k}=x_k 
\end{align*}
For polytopic $\mathbb{Z}$ and $\mathbb{X}_f$, this problem is a convex QP or SOCP if $\mathbb{W}$ is polytopic or ellipsoidal, respectively.
For all $x_k\in\mathbb{X}$, where $\mathbb{X}$ is the set of feasible states for $\mathcal{P}_1$, a receding control law is defined $u_k=v_{0|k}^*(x_k)$, where $(\cdot)^*$ denotes an optimal solution of $\mathcal{P}_1$. For all $x_0\in\mathbb{X}$, the closed-loop system
\begin{equation}\label{goulart_dynamics}
x_{k+1}=Ax_k+Bv_{0|k}^*(x_k)+Dw_k ,
\end{equation}
satisfies $x_k\in\mathbb{X}$ for all $k\in\mathbb{N}$~\cite{goulartandkerrigan2008}.
%
%
It is also shown in~\cite{goulartandkerrigan2008} that the system~\eqref{goulart_dynamics} is ISS with region of attraction $\mathbb{X}$.

\begin{prop}\label{asspropgoulart}
Assumptions \ref{distass}-\ref{terminallinearass} hold for the closed-loop system~(\ref{goulart_dynamics}) 
with $\Gamma = \Omega =\mathbb{X}_f$ and $\Phi = A+BK$.
\end{prop}

\bpf 
Assumption \ref{distass} holds due to the assumptions on $w_k$.
%
Assumptions \ref{terminalass} and \ref{terminallinearass} hold because 
$\Omega=\mathbb{X}_f\subseteq \mathbb{Z}$ is by assumption bounded and RPI for (\ref{lindyn}) under $u_k=Kx_k$, and since 
the solution of $\mathcal{P}_1$ is $v_{0|k}^\ast(x_k) = Kx_k$ for all $x_k\in\mathbb{X}_f$ due to the definition of $P$ and $K$.
%
%
Assumption~\ref{lyapunov} holds because (as proved in~\cite{goulartandkerrigan2008}) $J^\ast(x_k)$, the optimal value of the cost in $\mathcal{P}_1$, is an ISS-Lyapunov function for system \eqref{goulart_dynamics} satisfying the conditions of Theorem~\ref{lemiss}, thus guaranteeing that the system is ISS with region of attraction~$\mathbb{X}$.
\epf

Proposition~\ref{asspropgoulart} implies almost sure convergence to the minimal RPI set in~(\ref{eq:mRPIset}) and ensures bounds on average performance.

\begin{cor}
For all $x_0\in\mathbb{X}$, the closed-loop system \eqref{goulart_dynamics} satisfies \eqref{convlimit} and \eqref{perflimit} with $S=Q+K^\top R K$.
\end{cor}
\bpf This is a direct consequence of Corollary \ref{mainthm2} and Theorem \ref{mainthm3} since Assumptions \ref{distass}-\ref{terminallinearass} hold. \epf

Similar convergence results presented in~\cite{wangetal2008} required a redefinition of the cost and the control policy used in~\cite{goulartandkerrigan2008}. The results presented here apply to the original algorithm in~\cite{goulartandkerrigan2008} with the mild assumption that small disturbances have non-zero probability.

\subsection{Striped affine in the disturbance stochastic MPC}
\label{sec:kcandm}

The predicted control policy of~\cite{kouvaritakisetal2013} is again an affine function of future disturbance inputs. However, there are several differences with the formulation of Section~\ref{sec:gandk}: (i) a state feedback law with fixed gain is included in the predicted control policy; (ii) disturbance feedback gains are computed offline in order to reduce online computation; (iii) the disturbance feedback has a striped structure that extends over an infinite horizon. A consequence of (ii) is that state and control constraints can be enforced robustly by means of constraint tightening parameters computed offline, while (iii) has the effect of relaxing terminal constraints~\cite{kouvaritakisetal2013}.


The system is subject to mixed input-state hard and probabilistic constraints, defined for all $k\in\mathbb{N}$ by
\begin{subequations}
\label{koucons}
\begin{align}
    &(x_k,u_k)\in \mathbb{Z} \label{kouvaritakis_hcons},  \\
    &\mathbb{P}\{f^\top_j x_{k+1}+g^\top_j u_k\le h_j \}\ge p_j, \quad j\in\mathbb{N}_{[1,n_c]} \label{kouvaritakis_cons}
\end{align}
\end{subequations}
where $\mathbb{Z}$ is a convex compact polyhedral set containing the origin in its interior, $f_j\in \mathbb{R}^{n}$, $g_j\in \mathbb{R}^{n}$, $h_j\in\mathbb{R}$, $p_j\in(0, 1]$ and $n_c$ is the number of probabilistic constraints.

The predicted control policy has the structure
\begin{equation}\label{kouvaritakis_control_law}
u_{i|k}=
\begin{cases}
Kx_{i|k}+c_{i|k}+\sum_{j=1}^{i-1}L_{j}w_{i-j|k}, & i\in\mathbb{N}_{N-1} \\
Kx_{i|k}+\sum_{j=1}^{N-1}L_{j}w_{i-j|k}, & i \geq N
\end{cases}
\end{equation}
where ${\bf c}_k = \{c_{i|k}\}_{i\in\mathbb{N}_{N-1}}$ are decision variables at time $k$ and $K$ satisfies the algebraic Riccati equation $P=Q+A^\top PA-K^\top (R+B^\top PB)K$, $K=-(R+B^\top PB)^{-1}B^\top PA$. 

%
%

Using \eqref{lindyn} we note that constraint (\ref{kouvaritakis_cons}) can be imposed as 
\begin{equation}\label{tightened_prob}
    {f}^\top_j {x}_{k}+{g}^\top_j {u}_{k}+\gamma_j\le h_j,\quad j\in\mathbb{N}_{[1,n_c]}
\end{equation}
where $\tilde{f}^\top_j=f_j^\top A$, $\tilde{g}^\top_j=f_j^\top B+g_j^\top$ and $\gamma_j$ is a tightening parameter that accounts for the stochastic disturbance $w$. Specifically, $\gamma_j$ is computed using the PDF (or a collection of samples) of $w$ so that the satisfaction of the tightened constraint \eqref{tightened_prob} ensures satisfaction of \eqref{kouvaritakis_cons}. Constraint \eqref{tightened_prob} is equivalent to a polytopic constraint: $(x_k,u_k)\in{\mathbb S}=\{(x,u): {f}^\top_j {x}+{g}^\top_j \bar{u}+\gamma_j\le h_j, \, j\in\mathbb{N}_{[1,n_c]}\}$.

Let $\bar{x}_{i|k},\bar{u}_{i|k}$ be nominal state and input predictions satisfying $\bar{x}_{0|k}=x_k$, $\bar{x}_{i|k+1}=A\bar{x}_{i|k}+B\bar{u}_{i|k}$ and $\bar{u}_{i|k}=u_{i|k}$ (i.e.~assuming $w_{i|k}=0$ for all $i\in\mathbb{N}$). Then (\ref{koucons}) can be imposed with constraints on the nominal sequences
\begin{equation}
\label{tightened_constraints}
(\bar{x}_{i|k},\bar{u}_{i|k}) \in \tilde{\mathbb{S}}_i ,
\quad
(\bar{x}_{i|k},\bar{u}_{i|k}) \in \tilde{\mathbb{Z}}_i ,
\end{equation}
where $\tilde{\mathbb{S}}_i$ and $\tilde{\mathbb{Z}}_i$ are tightened versions of $\mathbb{S}$ and $\mathbb{Z}$. 
Here $\tilde{\mathbb{Z}}_i$ is computed by considering the worst-case of uncertain components of the predictions for the hard constraint $(x_{i|k},u_{i|k})\in\mathbb{Z}$, while $\tilde{\mathbb{S}}_i$ enforces
the probabilistic constraint \eqref{tightened_prob} and includes a worst-case tightening to ensure recursive feasibility. 


The disturbance feedback gains $L_{j}$, $j\in \mathbb{N}_{[1,N-1]}$, are computed sequentially  offline so as to minimize a the tightening of constraints \eqref{tightened_constraints}. Specifically, $L_1$ is first chosen so as to minimize the effect of the disturbance $w_{0|k}$ on the constraints \eqref{tightened_constraints} at prediction instant $i=2$, then $L_2$ is chosen so as to minimize the effect of $\{w_{0|k}, w_{1|k}\}$ on these constraints at prediction instant $i=3$, and so on until all $N-1$ gains have been chosen  (we refer the reader to~\cite{kouvaritakisetal2013} for further details).

%
The cost function is given by
\begin{equation}\label{kouvaritakis_cost}
J(x_k,{\bf c}_k)= \mathbb{E}\Bigl\{\sum_{i=0}^{\infty} (x_{i|k}^\top Qx_{i|k}+u_{i|k}^\top Ru_{i|k} -L_{ss} ) \Bigr\},
\end{equation}
where $Q,R\succ0$ and $L_{ss}=\lim_{k\rightarrow\infty}\mathbb{E} (x_{i|k}^\top Qx_{i|k}+u_{i|k}^\top Ru_{i|k})$ can be computed using the predicted control law for $i\geq N$ 
and the second moments of the disturbance input. It can be shown (e.g.~\cite[Thm.~4.2]{kouvaritakisetal2013}) that $J(x_k,{\bf c}_k)$ can be replaced, without changing the solution of the optimization problem, by
\begin{equation}\label{eq_eqcost}
V(x_k,{\bf c}_k) = x_k^\top P x_k+ {\bf c}_k^\top P_c {\bf c}_k
\end{equation}
where $P_c = I \otimes (R + B^\top P B)\succ 0$. 

The optimal control problem solved at the $k$th instant is therefore
\begin{align*}
\mathcal{P}_2: \ \ & \min_{{\bf c}_k} & & V(x_k,{\bf c}_k) \quad \text{s.t.} \quad  \ \forall  i\in\mathbb{N}_{N+N_2-1} \\
&&& \bar{u}_{i|k}=c_{i|k}+K\bar{x}_{i|k}\\
&&& (\bar{x}_{i|k},\bar{u}_{i|k}) \in \tilde{\mathbb{Z}}^{p}_i\cap \tilde{\mathbb{Z}}_i \\
&&& \bar{x}_{i+1|k}=A\bar{x}_{i|k}+B\bar{u}_{i|k} \\
&&& \bar{x}_{0|k}=x_k
\end{align*}
where $N_2$ is large enough to ensure that $(x_{i|k},u_{i|k})\in \tilde{\mathbb{S}}_i\cap\tilde{\mathbb{Z}}_i$ holds for all $i\geq N$. The solution ${\bf c}_k^\ast(x_k)$ defines the MPC law $u_k={c_{0|k}^*(x_k)+Kx_k}$
and the closed-loop dynamics are given, for $\Phi = A+BK$, by
\begin{equation}\label{dynamics_kouvaritakis}
x_{k+1}=\Phi x_k+Bc^*_{0|k}(x_k)+Dw_k .
\end{equation}
The set of states for which $\mathcal{P}_2$ is feasible, denoted by $\mathbb{X}$, is robustly invariant under the closed-loop dynamics~\cite{kouvaritakisetal2013}.

Asymptotically optimal performance is obtained if $u_k$ converges to the unconstrained optimal control law $u_k=Kx_k$. However, the bound
$\lim_{k\to\infty}\mathbb{E}\bigl\{x_{k}^\top Qx_{k}+u_{k}^\top Ru_{k}\bigr\} \leq L_{ss}$ is derived in~\cite[Thm~4.3]{kouvaritakisetal2013},
where $L_{ss} = l_{ss}+{\mathbb E}\{w^\top P_w w \}$ for ${P_w\succeq 0}$, and $l_{ss}=\lim_{k\rightarrow\infty}{\mathbb E}\{w^\top D^\top P D w\}$ is the asymptotic value of $\mathbb{E}\bigl\{x_{k}^\top Qx_{k}+u_{k}^\top Ru_{k}\bigr\}$ for~\eqref{lindyn} with $u_k=Kx_k$.
Thus, although~\cite{kouvaritakisetal2013} provides an asymptotic bound on closed-loop performance, this bound does not ensure convergence to the unconstrained optimal control law since $L_{ss} \geq l_{ss}$. 

From \eqref{eq_eqcost} it follows that the optimal solution is ${\bf c}_k^\ast(x_k) = 0$ whenever the constraints of $\mathcal{P}_2$ are inactive, and in particular we have ${\bf c}_k^\ast(0) = 0$. To ensure that the dynamics of (\ref{dynamics_kouvaritakis}) are linear on this system's limit set, we make the following assumption about the set $\mathbb{X}_{uc} = \{x\in\mathbb{X} : {c^\ast_{0|k}(x) = 0}\}$.

\begin{assum}\label{assterminalregimestriped}
The minimal RPI set (\ref{eq:mRPIset}) satisfies $\mathbb{X}_\infty \subseteq \mathbb{X}_{uc}$.
\end{assum}
It was proved in \cite{kouvaritakisetal2013} that Assumption \ref{assterminalregimestriped} guarantees the existence of an invariant set $\Omega$, such that $\mathbb{X}_\infty\subseteq \Omega\subseteq \mathbb{X}_{uc}$ and that if the state reaches $\Omega$, then it necessarily converges to $\mathbb{X}_\infty$. No guarantee is given in \cite{kouvaritakisetal2013} that the state will reach $\mathbb{X}_\infty$, but we can now apply the results of Sections~\ref{sec2} and~\ref{sec3} to the control policy of~\cite{kouvaritakisetal2013}.





\begin{thm}\label{cor1_kouvaritakis}
Let $V^\ast(x_k) = V\bigl(x_k,{\bf c}_k^\ast(x_k)\bigr)$, then $V^\ast(\cdot)$ is an ISS-Lyapunov function, and the closed-loop system \eqref{dynamics_kouvaritakis} is ISS with region of attraction $\mathbb{X}$.
\end{thm}

\bpf
First note that $\mathcal{P}_2$ is a convex QP since ${P,P_c\succ 0}$. Therefore $V^\ast(\cdot)$ is strictly convex and continuous and the optimizer ${\bf c}_k^\ast(\cdot)$ is also continuous~\cite[Thm.~4]{bemporad02}. It follows that (i) $V^\ast(\cdot)$ is Lipschitz continuous on $\mathbb{X}$, (ii) the function $f(x_k,w_k) = \Phi x_k+Bc_{0|k}^*(x_k)+Dw_k$ defining the closed-loop dynamics in \eqref{dynamics_kouvaritakis} is continuous on $x_k$, and (iii) condition (\ref{lyapkappa}) holds since $V^\ast(\cdot)$ is positive definite. 
Next, from~\eqref{dynamics_kouvaritakis}, \eqref{eq_eqcost} and $P-\Phi^\top P \Phi = Q + K^\top R K$ we have $V(f(x_k,0))-V(x_k) \leq -(x_k^\top Qx_k + u^{\top}_kRu_k)\le-x_k^\top Qx_k$. Then, there exists a ${\mathcal K}_\infty$ function $\alpha_3(\cdot)$ such that $V\left(f(x,0)\right)-V(x)\le -\alpha_3(\norm{x})$ holds since $Q\succ 0$. Applying \cite[Lem.~22]{goulartetal2006}, these conditions imply that \eqref{dynamics_kouvaritakis} is ISS with region of attraction $\mathbb{X}$. 
\epf

\begin{prop}\label{asspropkovaritakis}
Assumptions \ref{distass}-\ref{terminallinearass} hold for the closed-loop system (\ref{dynamics_kouvaritakis}) with $\Gamma = \Omega  \subseteq \mathbb{X}_{uc}$
if Assumption~\ref{assterminalregimestriped} holds.
\end{prop}

\bpf Assumption \ref{distass} due to the assumptions on $w_k$.
%
Assumptions~\ref{terminalass} and~\ref{terminallinearass} hold because $\mathbb{X}_\infty$ is RPI, $\mathbb{Z}$ is a bounded set, and ${\bf c}_k^\ast(x_k) = 0$ for all $x_k\in\mathbb{X}_\infty$ under Assumption~\ref{assterminalregimestriped}. Assumption~\ref{lyapunov} holds by Theorem~\ref{cor1_kouvaritakis}.
\epf

This allows us to conclude the following convergence results for the state and limit average performance.

\begin{cor}
For all $x_0\in\mathbb{X}$, the closed-loop system \eqref{dynamics_kouvaritakis} satisfies \eqref{convlimit} and \eqref{perflimit} with $S = Q + K^\top R K$ 
\end{cor}
\bpf 
The bounds in \eqref{convlimit} and \eqref{perflimit} follow from Corollary~\ref{mainthm2} and Theorem~\ref{mainthm3}  since Assumptions \ref{distass}-\ref{terminallinearass} hold. \epf

\subsection{Nonlinear stochastic MPC based on constraint-tightening}\label{nlsmpc}
This section considers the convergence properties of a nonlinear system with the stochastic MPC algorithm~\cite{santosetal2019}. 
Assuming that arbitrarily small disturbances have a non-vanishing probability, we use the results of~\cite{santosetal2019} and Section \ref{sec3} to show that the closed-loop system converges with probability~1 to the minimal RPI set associated with the MPC law.

We consider the system with state $x\in{\mathbb R}^n$, disturbance input $w\in\mathbb{W}\subset\mathbb{ R}^{n_u}$ and control input $u\in{\mathbb R}^n$ governed by
\begin{equation}\label{nldynamics}
    x_{k+1}=f(x_k,u_k)+w_k .
\end{equation}
The function $f(\cdot,\cdot)$ satisfies $f(0,0)=0$ and is assumed to be uniformly continuous in its arguments for any feasible pair $(x,u)$. 
%
%
The system is subject to chance constraints on its state and mixed input-state hard constraints of the form
\begin{subequations}
\label{conssantos}
\begin{align}
    &(x_k,u_k)\in \mathbb{Z} \label{santoshardc}, \\ 
    \label{santosprobc}
    &\mathbb{P}\{g^\top_j x_{k+1}\le h_j \}\ge p_j, \quad j\in\mathbb{N}_{[1,n_c]}
\end{align}
\end{subequations}
with $g\in \mathbb{R}^{n}$, $h\in\mathbb{R}$,  $p\in(0, 1]$. Here $n_c$ is the number of probabilistic constraints and for simplicity we assume that the projection $\{x : \exists u, \, (x,u) \in \mathbb{Z}\}$ is bounded.
The predicted control law is parameterized by decision variables $v_{i|k}$ so that $u_{i|k}=\pi(x_{i|k},v_{i|k})$, where $\pi(\cdot,\cdot)$ is assumed to be continuous on the feasible domain. Defining $f_\pi(x,v)=f(x,\pi(x,v))$, the predicted states  therefore evolve according to
\begin{equation}\label{nlcldynamics}
    x_{i+1|k}=f_\pi(x_{i|k},v_{i|k})+w_{i|k}.
\end{equation}
The cost and constraints in the formulation of the optimal control problem are defined over the nominal state predictions $\bar{x}_{i|k}$, obtained from \eqref{nlcldynamics} but assuming $w_{i|k}=0$ for all $i$.

Noting that \eqref{santosprobc} can be imposed as a polyhedral constraint $x_{k+1}\in{\mathbb S}$, constraints \eqref{conssantos} are imposed in the optimal control problem for a prediction horizon $N$ by
\[
(x_{i|k},v_{i|k})\in \tilde{{\mathbb{Z}}}_i, 
\quad 
x_{i+1|k}\in\tilde{\mathbb S}_i,
\]
where $\tilde{\mathbb{Z}}_i$ and $\tilde{\mathbb S}_i$ are tightened versions of $\mathbb{Z}$ and ${\mathbb S}$. These are computed by considering worst-case disturbances for the hard constraints (\ref{conssantos}a) and a combination of worst-case and stochastic disturbances to ensure recursive feasibility of the probabilistic constraints (\ref{conssantos}b). The construction of the required tightened sets relies on the uniform continuity of $f(x,y)$ for any pair $(x,y)$ and can be performed offline.  For further details we refer the reader to \cite{santosetal2019}.


The stochastic MPC formulation assumes a terminal control law $v_t(\cdot)$ (so that $u_{i|k} = \pi\bigl(x_{i|k},v_t(x_{i|k})\bigr)$ for $i\geq N$), a terminal constraint ${\mathbb X}_f$ and a terminal cost $V_f(\cdot)$ that satisfy the  following conditions. (i) $\pi\bigl(0,v_t(0)\bigl)=0$.
(ii) ${\mathbb X}_f$ is RPI, with 
$f_\pi(x,v_t(x)) + w \in{\mathbb X}_f$ for all $(x,w)\in{\mathbb X}_f\times\mathbb W$, 
and ${\mathbb X}_f\subseteq\tilde{\mathbb S}_N\cap\{ x:(x,v_t(x))\in\tilde{\mathbb Z}_N\}$. (iii) $V_f(x)$ is a positive definite function satisfying, for all $x,x_1,x_2\in{\mathbb X}_f$: $\alpha_{1f}(\|x\|)\le V_f(x)\le \alpha_{2f}(\|x\|)$, $V_f(x,v_t(x))-V_f(x)\le-L_\pi(x,v_t(x))$ and $|V_f(x_1)-V_f(x_2)|\le\delta(\|x_1-x_2\|)$, where $\alpha_{1f}(\cdot)$, $\alpha_{2f}(\cdot)$ and $\delta(\cdot)$ are $\mathcal{K}$-functions.


The cost function is given by
\begin{equation}\label{santos_cost}
J(x_k,{\bf v}_k)=V_f(\bar{x}_{N|k}) +  \sum_{i=0}^{N-1} L_\pi(\bar{x}_{i|k},v_{i|k} ),
\end{equation}
where ${\bf v}_k = \{v_{i|k}\}_{i\in\mathbb{N}_{N-1}}$ and $L_\pi(x,v)$ is a positive definite and uniformly continuous function for all feasible $(x,v)$.
The optimal control problem solved at the $k$th instant is
\begin{align*}
\mathcal{P}_3: \ \ & \min_{{\bf v}_k} & & J(x_k,{\bf v}_k) \quad \text{s.t.} \quad \forall i\in\mathbb{N}_{N-1} \\
&&& \bar{x}_{i+1|k}\in \tilde{\mathbb S}_i \\
&&& (\bar{x}_{i|k},v_{i|k}) \in \tilde{\mathbb{Z}}_i \\
&&& \bar{x}_{i+1|k}=f_\pi(\bar{x}_{i|k},v_{i|k}) \\
&&& \bar{x}_{0|k}=x_k, \; \bar{x}_{N|k}\in {\mathbb X}_f.
\end{align*}
Let $\mathbb{X}$ be the set of all states $x_k$ such that ${\mathcal P}_3$ is feasible. Then for all $x_k\in\mathbb{X}$ the solution $v^*_{0|k}(x_k)$ defines the MPC law $u_k=\pi(x_k,v_{0|k}^*(x_k))$
and the closed-loop system is given by 
\begin{equation}\label{santosclosedloop}
x_{k+1}=f_\pi(x_k,v_{0|k}^*(x_k))+w_k.
\end{equation}
It can be shown~\cite[Lemma 3.1]{santosetal2019} that $x_{k+1}\in\mathbb{X}$ if $x_{k}\in\mathbb{X}$, so the constraints \eqref{conssantos} are satisfied  for all $k\in\mathbb{N}$ if $x_0\in\mathbb{X}$, and furthermore the system~\eqref{santosclosedloop} is ISS~\cite[Appendix B]{santosetal2019}. 
%
%
%

To find a tight limit set, let ${\mathbb X}_\infty$ be the minimal RPI set of \eqref{santosclosedloop}. The properties of this closed-loop system that allow the results of Section \ref{sec3} to be applied are summarized as follows.


\begin{prop}
Assumptions \ref{distass}-\ref{lyapunov} hold for the closed-loop system \eqref{santosclosedloop} and the set
$\Omega =\mathbb{X}_\infty$.
\end{prop}
\bpf
Assumption \ref{distass} holds due to the assumptions on $w_k$, Assumption \ref{terminalass}  holds because
$\mathbb{X}_\infty$ is bounded due to $\mathbb{X}_\infty\subseteq\mathbb{Z}$ and
Assumption~\ref{lyapunov} follows because \eqref{santosclosedloop} is ISS. \epf 

Corollary \ref{mainthm2} therefore implies that the state of the closed-loop system (\ref{santosclosedloop}) converges to ${\mathbb X}_\infty$ with probability~1.

%


\subsection{Discussion}

\textit{On the use of the proposed analysis for MPC of linear and nonlinear systems:} We first highlight that even for linear systems under a stochastic MPC law it is necessary, in order to use the results of Section \ref{sec3} to analyse convergence, that the analysis is posed for nonlinear systems. In fact, even though system \eqref{lindyn} is linear, the closed-loop dynamics under the stochastic MPC laws in Sections~\ref{sec4}.A and \ref{sec4}.B are in general (except on the MPC terminal set) nonlinear. We note that stability analyses based on ISS properties have been developed for several nonlinear stochastic MPC formulations in the literature~\cite{sehrbitmead2017,santosetal2019,lorenzenetal2019}. This suggests that analysis discussed here is a generally useful tool for stochastic MPC.



\textit{On the MPC law in the limit set:} We have shown that the strategies of Sections \ref{sec4}.A and \ref{sec4}.B recover their respective optimal control laws in the limit set, and their limit average performance is optimal for the MPC cost function.
%
These results are obtained even though the terminal predicted control law of the MPC formulation in Section \ref{sec4}.B incorporates sub-optimal disturbance compensation terms and is thus sub-optimal with respect to the MPC cost. In fact, the MPC law in the limit set is not necessarily the same as the terminal control law of the MPC formulation. 
This highlights that in order to design a stochastic MPC strategy with optimal limit average performance it is not necessary that the terminal control law is optimal. Instead, we require that MPC law equals the optimal control law on an RPI set that contains the origin in its interior.

The procedure consists of checking Assumptions~\ref{distass}-\ref{lyapunov} and invoking Corollary \ref{mainthm2} to conclude that the state converges to the minimal RPI set for the closed-loop system. Furthermore, if the dynamics of the closed-loop system are linear within some RPI set then a tight limit set can be determined explicitly (e.g.~\cite{blanchini08}).

\section{Conclusions}

This paper extends and generalizes methods for analysing the convergence for disturbed nonlinear systems, which can be applied to stochastic MPC formulations. 
We define a set of conditions 
on stochastic additive disturbances (Assumption~\ref{distass}); 
on the existence of an invariant set (Assumption~\ref{terminalass}) and an ISS-Lyapunov function (Assumption~\ref{lyapunov}) for the closed-loop system.
%
We show that a nonlinear stochastic system satisfying these conditions converges almost surely to a limit set defined by the minimal RPI set of the system. For the case where the dynamics are linear in the limit set, we show that the asymptotic average performance is tightly bounded by the performance associated with those linear dynamics. 
The results are obtained using the ISS property of the system, but the limits directly implied by the ISS Lyapunov inequality would yield worse asymptotic bounds. 
These conditions are commonly met by stochastic MPC strategies, and we illustrate the use of the convergence analysis by applying it to three existing formulations of stochastic MPC. 
In each of these applications our analysis allows for 
improved ultimate bounds on state and performance.

\bibliographystyle{ieeetr}
\bibliography{terminalconv_v14}

\end{document}